\documentclass[11pt,a4paper]{amsart}

\usepackage{geometry}
\usepackage[utf8]{inputenc}
\usepackage{amsmath,amssymb,amsthm}
\usepackage[pagebackref=false]{hyperref}
\hypersetup{
colorlinks = true,
pdfstartview = FitH,
pdfpagemode = UseNone,
pdfauthor    = {Giacomo Cherubini, Alessandro Fazzari, Andrew Granville, Vitezslav Kala, Pavlo Yatsyna},
pdftitle     = {},
pdfkeywords  = {},
pdfcreator   = {Pdflatex},
linktoc = none,
}

\numberwithin{equation}{section}

\def\eps{\varepsilon}

\def\QQ{\mathbb{Q}}

\theoremstyle{plain}
\newtheorem{theorem}{Theorem}
\newtheorem{lemma}[theorem]{Lemma}

\theoremstyle{definition}

\title{Consecutive real quadratic fields with large class numbers}
\author[Cherubini]{Giacomo Cherubini}
\author[Fazzari]{Alessandro Fazzari}
\author[Granville]{Andrew Granville}
\author[Kala]{V\'\i t\v ezslav Kala}
\author[Yatsyna]{Pavlo Yatsyna}
\date{\today}
\address{
         Charles University,
         Faculty of Mathematics and Physics,
         Department of Algebra,
         Sokolov\-sk\'a 83, 18600 Praha~8,
         Czech Republic
        }
\email{
    cherubini@karlin.mff.cuni.cz\\ 
    fazzari@karlin.mff.cuni.cz\\
    kala@karlin.mff.cuni.cz\\ 
    yatsyna@karlin.mff.cuni.cz
    }

\address{
         D\'epartment de Math\'ematiques et Statistique,
         Universit\'e de Montr\'eal, CP 6128
         succ Centre-Ville, Montr\'eal, QC H3C 3J7,
         Canada
        }

\email{andrew.granville@umontreal.ca}

\thanks{This work was supported by Czech Science Foundation GACR, grant 21-00420M (GC, AF, VK), the project PRIMUS/20/SCI/002 from Charles University (GC, VK, PY), and Charles University Research Centre program UNCE/SCI/022 (GC, VK). AG was supported by grants from NSERC (Canada)}
\keywords{Class number, continued fraction, real quadratic field, special $L$-values}
\subjclass[2010]{Primary 11R29; Secondary 11R11, 11A55, 11M20.}
\begin{document}

%%  ABSTRACT AND TITLE %%
\begin{abstract}
For a given positive integer $k$, we prove that there are at least $x^{1/2-o(1)}$ integers $d\leq x$
such that the real quadratic fields
\(\QQ(\sqrt{d+1}),\dots,\QQ(\sqrt{d+k})\)
have class numbers essentially as large as possible.
\end{abstract}

\maketitle

%% SECTION : INTRODUCTION %%
\section{Introduction}

The famous class number one problem for real quadratic fields
states that infinitely many of them should have class number one.
As of today, this is still widely open.
Nevertheless, quadratic fields with class number one
were classified in the imaginary case \cite{heegner}
and for specific families (such as Yokoi's and Chowla's)
in the real case~\cite{biro:yokoi,biro:chowla,BG,BL,Lo,Mo}.
In contrast, it is known that
the class number can be arbitrarily large (see, e.g.~\cite{abc,montgomery}),
and a great deal of work was done to understand its behaviour.
Explicit results about the distribution of class numbers over families
of real quadratic fields were obtained in \cite{dahl-kala,dahl-lamzouri}.

Lamzouri \cite{lamzouri} gave explicit constants in the work of Montgomery and Weinberger \cite{montgomery}, showing that 
\[
h_{\QQ(\sqrt{d}\,)} \geq (2e^\gamma+o(1)) \frac{\sqrt{d}}{\log d}\log\log d, 
\]
for $x^{1/2+o(1)}$ real quadratic fields $\QQ(\sqrt{d}\,)$ with $d\leq x$.
We believe this is  as large as possible; certainly class numbers cannot be much larger since, assuming the Generalized Riemann Hypothesis, we know that \cite{littlewood}
\begin{equation}\label{eq:grh}
	h_{\QQ(\sqrt{d}\,)} \leq (4e^\gamma+o(1)) \frac{\sqrt{d}}{\log d}\log\log d.
\end{equation}

The behaviour of the class number is rather erratic in general:
as one visits consecutive integers $d,d+1,d+2,\dots$, the associated real quadratic
fields can have very different class numbers.

Below we present a heuristic which suggests that there are only
$x^{1/2+o(1)}$ values of  $d\leq x$ for which $h_{\QQ(\sqrt{d}\,)}\gg \frac{\sqrt{d}}{\log d}\log\log d$; one might guess that these
should be roughly randomly distributed so it is unlikely that there are sequences of consecutive $d$-values
with \emph{large} class numbers. Nevertheless, our main theorem shows the existence of such $k$-tuples:

\begin{theorem}\label{TG}
Fix an integer $k\geq 1$. There are $\geq x^{1/2-o(1)}$ integers $d\leq x$ such that
\[
h_{\QQ(\sqrt{d+j}\,)} \gg_k \frac{\sqrt{d}}{\log d}\log\log d,\quad \forall\;j=1,\dots,k.
\]
\end{theorem}

A key ingredient in the proof of Theorem \ref{TG} is that if $d$ is suitable chosen,
then for all $j$,
the fundamental units in the fields $\QQ(\sqrt{d+j})$ are bounded by $d$ times a constant.
At the same time we can show, by a combination of the Chinese Remainder Theorem
and a sieving argument, that for many values of $d$
the discriminant $D_j$ of such fields and the associated value $L(1,\chi_{D_j})$ are essentially as large as $d$ and $\log\log d$, respectively
(here $\chi_{D_j}$ is the primitive quadratic character modulo $D_j$).
Using these two inputs in Dirichlet's class number formula we obtain our theorem.

In view of the conditional upper bound \eqref{eq:grh}, Theorem \ref{TG} is the best possible up to the implied constant. Indeed, from  GRH it follows that $L(1,\chi_{D})\ll \log\log D$ for every fundamental discriminant $D$
(see, e.g.~\cite{littlewood}), which is used to prove \eqref{eq:grh}.  We obtain the reverse inequality $L(1,\chi_D)\gg_k \log\log D$ unconditionally.

\medskip

\subsection*{Heuristic} Finding many consecutive integers with small fundamental units contradicts the ``obvious'' heuristic. Given positive integers $D$
and $X$, let us consider the set
\begin{equation}\label{ale1} \{d\asymp D : X<\varepsilon_d<2X\}. \end{equation}
To estimate the cardinality of this set, we write $\varepsilon_d=u+v\sqrt{d}$; the condition $X<\varepsilon_d<2X$ yields $u\asymp X$ and $v\asymp X/\sqrt{d}$. 
As $\varepsilon_d$ solves the Pell equation, we have \begin{equation}\label{ale2} dv^2\pm 1 = u^2. \end{equation}
For the left hand side, we have about $X/\sqrt{d}$ 
values for $v$ such that $dv^2\asymp X^2$, and the probability that a number of size $X^2$ is a square is $1/X$.
Therefore, the probability that for a fixed  $d$ there is some $v$ satisfying \eqref{ale2}  is approximately $1/\sqrt{d}$. Hence, we expect that the set in \eqref{ale1} has size roughly $D^{1/2+o(1)}$.

\medskip

Note that we did not explicitly place any restriction on the sizes of $D$ or $X$ in the very rough heuristic above. However, a reasonable range is, e.g., $\sqrt{D}/2<X<e^{O(\sqrt{D})}$, as
if we assume GRH, then 
$\log\eps_D=\frac{\sqrt{D}L(1,\chi_{D})} {h_{\QQ(\sqrt{D})}}\leq \sqrt{D}L(1,\chi_{D}) \ll \sqrt{D} \log\log D$,
and so $\sqrt{D}/2< \eps_D< (\log D)^{O( \sqrt{D} )}$.

\medskip

In particular, we expect that $\log \varepsilon_d\ll \log d$ for $\ll x^{1/2+o(1)}$ integers $d\leq x$. A naive guess is that these are ``randomly distributed'', but if we have $x^{1/2+o(1)}$ randomly selected integers $\leq x$, we  expect to have very few tuples of consecutive integers amongst them, and certainly not $\geq x^{1/2-o(1)}$ $k$-tuples.

\medskip

Let us also remark that several recent papers have studied pairs \cite{iizuka}, triples \cite{cm},
or in general $k$-tuples \cite{ikn} of quadratic fields in relation to the divisibility
property of the class number. Regarding this, we mention the work of Hoque \cite{hoque},
who showed that there is an infinite family of quintuples of imaginary quadratic fields
\[
\QQ(\sqrt{d}),\QQ(\sqrt{d+1}),\QQ(\sqrt{d+4}),\QQ(\sqrt{d+36}),\QQ(\sqrt{d+100})
\]
whose class numbers are all divisible by a given integer.
Theorem \ref{TG} should be compared to Iizuka's conjecture \cite{iizuka}
that for any prime $l$ and any positive integer $n$, there should be an
infinite family of quadratic fields
\[
\QQ(\sqrt{d}),\QQ(\sqrt{d+1}),\dots,\QQ(\sqrt{d+k})
\]
whose class numbers are all divisible by $l$.

\medskip

Finally, in Section \ref{S3} we address the question of the uniformity in $k$ in Theorem \ref{TG}
and prove that the implied constant can be taken of the form $Ck^{-5/2}$ for some absolute constant $C>0$.
The exponent of $k$ can be improved if we do a geometric average of the class numbers:
in Section \ref{S3} we also prove that
\begin{equation}\label{1812:eq002}
	\bigg(\prod_{j=1}^k h_{\QQ(\sqrt{d+j})}\bigg)^{1/k} \geq \frac{C'\sqrt{d}\log\log d}{k^2\log d}.
\end{equation}

\medskip

Recall that
for non-negative functions $f(x),g(x)$ we write $f\ll g$ (or $f=O(g)$) if there exists
a positive constant $C$ such that $f(x)\leq Cg(x)$ for all $x$ sufficiently large, and $f=o(g)$ if $\lim f(x)/g(x)=0$.
Further, we will use the notations $f\ll_P g, f=O_P(g)$ to stress that the constant $C$
depends on the parameter(s) $P$.

\section*{Acknowledgments}

We wish to thank Youness Lamzouri for several helpful suggestions.

\section{Proof of Theorem \ref{TG}}\label{S2}

As anticipated in the introduction, we want to apply Dirichlet's class number formula.
Therefore, our goal is to find consecutive integers for which we can control
the size of the regulator and of the discriminant in the associated quadratic fields,
as well as the associated $L$-functions.

\medskip

To this purpose, we place ourselves in a one-parameter family:
given $k\geq 1$, define $P=P(k)=\mathrm{lcm}[1,\dots,k]$ and take
\begin{equation}\label{def:delta}
\Delta(m)=(mP)^2.
\end{equation}
If we set $P_j=P/j$ for $1\leq j\leq k$, then we can write
\begin{equation}\label{def:Deltamj}
\Delta(m)+j = (mP_j)^2 j^2+j.
\end{equation}
The choice of $\Delta(m)$ as in \eqref{def:delta} is motivated by the fact 
that,
for any positive integer of the form $N=A^2j^2+j$, we have the identity
$(2A^2j+1)^2-N(2A)^2=1$, so the element
\[
\eps = 2A^2j+1+2A\sqrt{N}
\]
is a unit greater than one in $\QQ(\sqrt{N})$.
Since $\eps$ must be a power of the fundamental unit $\eps_N$,
it follows that the regulator of the field is
\begin{equation}\label{bound:regulator}
\log \eps_N \leq \log \eps \leq  \log 5N.
\end{equation}
In particular, by \eqref{def:Deltamj}, this is true
for $N=\Delta(m)+j$, for all $j=1,\dots,k$.

\medskip

It is worth noting that the existence of a small unit,
and thus the bound \eqref{bound:regulator} on the regulator,
reflects the fact that the continued fraction of $\sqrt{\Delta(m)+j}$
is particularly simple. Indeed, for any integer $N=A^2j^2+j$, we have
\[
\sqrt{N} = [Aj, \overline{A, 2Aj}].
\]
By the theory of continued fractions (see, e.g.~\cite[Corollary 3.3.2]{aa}),
we know that a suitable convergent provides a solution to the Pell equation $x^2-Ny^2=1$.
In turn, this gives a unit $\eps$ in $\QQ(\sqrt{N})$ of polynomial size
and therefore the regulator is at most a constant times $\log N$.
We note also that integers of the form $A^2j^2+j$ are a special case of those
studied by Schinzel \cite{schinzel} and Friesen \cite{friesen}
(actually, their results led us to our choice of $\Delta(m)$).

\medskip

Next, regarding the $L$-functions, we know by \cite[Proposition 2.2]{granville-sound} that
for all but $Q^{2/7+o(1)}$ of the primitive characters $\chi$ ($\mathrm{mod}\;q$) with $q\leq Q$
we have
\begin{equation}\label{eq:L}
L(1,\chi) = (1+o(1)) \prod_{p\leq (\log Q)^7} \left(1-\frac{\chi(p)}{p}\right)^{-1}.
\end{equation}
Therefore, in order to make the above large,
we would like to have $\chi(p)=1$ for all primes up to $(\log Q)^7$.
In reality, for any $\epsilon\in(0,1)$,
we can focus on primes up to $(\log Q)^\epsilon$, since
the product over primes
in the interval $J=[(\log Q)^\epsilon,(\log Q)^7]$
can be bounded from below,
using Merten's theorem, by
\[
\prod_{p\in J} \left(1-\frac{\chi(p)}{p}\right)^{-1}
\geq
\exp\bigg(-\sum_{p\in J} \frac{1}{p} + O(1)\bigg)
\gg
\exp(O_\epsilon(1)),
\]
where the implied constant does not depend on $Q$.
From now on, we assume thus that we have fixed a value for $\epsilon\in(0,1)$.

\medskip

We wish to show, using \eqref{eq:L}, that for many values of $m$ the $L$-values are large;
at the same time, we want the discriminant $D_j$ of $\QQ(\sqrt{\Delta(m)+j})$
to be essentially of the same size as $\Delta(m)+j$.

\medskip

Let us consider integers $m$ such that $\Delta(m)\leq x$, i.e. $m\leq M:=\sqrt{x}/P$.
Also, we define
\begin{equation}\label{def:q}
q=\prod_{k<p\leq (\log x)^\epsilon} p.
\end{equation}

Let us now start by showing that
the characters defining the $L$-values can be assumed to be all equal to one,
provided $m$ lies in a suitable residue class modulo $q$.

\medskip

\begin{lemma}\label{lemma1}
Let $k\geq 1$ and let $D_j$ and $q$ be as above.
Then, there exists $m_0\pmod{q}$ such that, if $m\equiv m_0\pmod{q}$, we have
\[
\chi_{D_j}(p) = 1, \qquad \forall\; 5^k<p\leq (\log x)^{\epsilon}.
\]
Moreover, we can take $m_0\equiv 0\pmod{p}$ for $k<p\leq 5^k$.
\end{lemma}

\begin{proof}
The strategy is to start with
an arithmetic progression $m_0$ $(\mathrm{mod}\;p)$ for every prime
$k<p\leq (\log x)^\epsilon$
and then glue everything by the Chinese Remainder Theorem,
which yields an arithmetic progression$\pmod{q}$.

To construct $m_0$, we first impose $m_0\equiv 0$ $(\mathrm{mod}\;p)$ for primes $p$ with
$k<p\leq 5^k$ and then we select $m_0$ $(\mathrm{mod}\;p)$ so that
\begin{equation}\label{0512:eq001}
\bigg( \frac { (m_0P)^2+1}{p}\bigg) = \bigg( \frac { (m_0P)^2+2}{p}\bigg) = \cdots = \bigg( \frac { (m_0P)^2+k}{p}\bigg) = 1
\end{equation}
for all primes $p$ with $5^k<p\leq (\log x)^\epsilon$.
We need to show that at least one such $m_0$ $(\mathrm{mod}\;p)$ exists.

Let $n\equiv (m_0P)^2 \pmod p$. Note that $P$ is invertible$\pmod{p}$
for primes $p>k$.
Therefore, there are twice as many non-zero classes $m_0\pmod p$
solving \eqref{0512:eq001} as~the number of $n \pmod p$ with 
\begin{equation}\label{0612:eq001}
\bigg(\frac{n}{p}\bigg) = \bigg(\frac{n+1}{p}\bigg) =  \cdots = \bigg(\frac{n+k}{p}\bigg) = 1,
\end{equation}
where the first Legendre symbol accounts for the fact that $n$ is a square$\pmod{p}$.
The number of solutions to \eqref{0612:eq001} is
\begin{equation}\label{0512:eq002}
\frac 1{2^{k+1}} \sum_{n=1}^{p-k-1} \prod_{j=0}^k \bigg( 1 + \bigg( \frac { n+j}{p}\bigg)\bigg)
=
\frac 1{2^{k+1}} \sum_{n=1}^p \prod_{j=0}^k \bigg( 1 + \bigg( \frac { n+j}{p}\bigg)\bigg) +O(k).
\end{equation}
We multiply this product out and then sum over $n$; the main term comes from all the $1$s, and gives $p/2^{k+1}$.
The other terms are all of the form
$2^{-k-1} \sum_n (\frac{f(n)}{p})$, with $f(x)$
being a polynomial of degree at least one and at most $k$.
Weil (\cite{weil}, see also \cite[Theorem 11.23]{ik}) proved
that each of these sums is $\leq k \sqrt{p}$, and we have
no more than $2^{k+1}$ such sums.
Hence, \eqref{0512:eq002} equals
\[
\frac{p}{2^{k+1}} +O(k\sqrt{p}) > \frac{p}{2^{k+2}} ,
\]
since $p>5^k$.
Therefore, there is at least one (and in fact many) $m_0\pmod{p}$ solving \eqref{0512:eq001}
for each such prime.
Lifting by means of the Chinese Remainder Theorem,
we obtain at least one residue class$\pmod{q}$, as desired.
\end{proof}

In view of the above lemma, we will restrict ourselves to integers $m$
in a specific arithmetic progression$\pmod{q}$.
For these values of $m$, we want to control the size of $D_j$.

\begin{lemma}\label{lemma2}
Let $q$ be as in \eqref{def:q} and let $m_0\pmod{q}$ be
as given by Lemma~\ref{lemma1}.
For $M^{1-o(1)}$ integers $m\leq M$ with $m\equiv m_0$, we have
\[
D_j \gg_k \Delta(m) \qquad \forall\;j=1,\dots,k.
\]
\end{lemma}

\begin{proof}
We start by setting some notation:
let $S_j^2$ be the largest square dividing~$\Delta(m)+j$.
Then, recalling that the discriminant of $\QQ(\sqrt{N})$ is equal to the squarefree part of $N$
or four times that, we have $S_j^2D_j=\Delta(m)+j$ or $4(\Delta(m)+j)$.
Furthermore, let $s_j^2$ be the largest square dividing $j$.
In particular, we have $s_j|S_j$, and so we can write $S_j=r_js_j$
for some positive integer $r_j$.

In what follows we will show that,
for $M^{1-o(1)}$ integers $m\equiv m_0\pmod{q}$, we have $r_j=1$.
Since $s_j\leq \sqrt{k}$, this will give
\begin{equation}\label{def:sj}
D_j \gg \frac{\Delta(m)+j}{s_j^2r_j^2} \gg_k \Delta(m).
\end{equation}

By definition (see \eqref{def:Deltamj}), $\Delta(m)+j$ is divisible by $j$.
If we define the polynomial
\[
F_{k,j}(t) = jP_j^2t^2 +1,
\]
then we can write $\Delta(m)+j=j F_{k,j}(m)$,
and it suffices to study $F_{k,j}(m)$.

First, note that $F_{k,j}(m)$ is not divisible by any prime $p\leq 5^k$,
because
$p$ divides $jP_j$ for $p\leq k$ and
$m\equiv m_0\equiv 0\pmod{p}$ for $k<p\leq 5^k$.
Also, if $m\equiv m_0\pmod{q}$, then by construction we have $\big(\frac{\Delta(m)+j}{p}\big)=1$
for all $p|q$, and so in particular $p$ does not divide $F_{k,j}(m)$.
Therefore, if $m\equiv m_0\pmod q$,
the prime factors of $F_{k,j}(m)$ are all larger than $(\log x)^\epsilon$.

For primes $p>(\log x)^\epsilon$, it could still be that $\Delta(m)+j$
is divisible by $p^2$ (which would imply $r_j>1$). However,
by a sieving argument we can show that this does not happen too often.

Denoting $z=q^2(\log x)^{4k}$, we claim that for $\gg_k M^{1-o(1)}$
integers $m\leq M$ in the arithmetic progression $m_0\pmod{q}$ we have:

$(i)$ $p$ does not divide any $F_{k,j}(m)$, for all primes $p$ with  $(\log x)^\epsilon <p\leq z$ and all $j=1,\dots,k$. 

$(ii)$ $p^2$ does not divide any $F_{k,j}(m)$, for all primes $p$ with $z <p\leq 2MP/z^{1/2}$ and all $j=1,\dots,k$.

Regarding $(i)$, we use the small sieve.
There are at most two congruence classes for $m$ such that
$F_{k,j}(m)\equiv 0 \pmod p$, and  $k$ polynomials.
Therefore, we need to sieve $g(p)$ residue classes, where $g(p)\leq 2k$.
The number of integers left unsieved is
\begin{equation}\label{0612:eq012}
\gg
\frac Mq \prod_{(\log x)^\epsilon <p\leq z}  \bigg( 1 - \frac {g(p)}{p}\bigg)  
\geq
\frac Mq \prod_{(\log x)^\epsilon <p\leq z}  \bigg( 1 - \frac {2k}{p}\bigg)
\gg_k
\frac M{q(\log z)^{2k}},
\end{equation}
As for $(ii)$, observe that for each $p$
there are at most $2k$ congruence classes $m\pmod {p^2}$ for which $p^2$ divides $F_{k,j}(m)$
for some $j$. Therefore, the number of such $m$ is bounded by
\[
2k\bigg(1 + \frac M{p^2q}\bigg),
\]
for each given prime $p$.
The total number of $m$ removed in this way is thus bounded by
\begin{equation}\label{0612:eq013}
\ll_k
\sum_{z <p\leq 2MP/z^{1/2}} 2k\bigg(1 + \frac{M}{p^2q}\bigg)
\ll_k
\frac{M}{z^{1/2}\log M}  + \frac{M}{q z\log z}.
\end{equation}
By our choice of $z$, we deduce that \eqref{0612:eq012} is larger than \eqref{0612:eq013}.
Moreover, since $q=x^{o(1)}=M^{o(1)}$, we conclude that there are at least $M^{1-o(1)}$
integers for which $(i)$ and $(ii)$ hold.

In order to conclude the proof of the lemma, assume that
$p^2|F_{k,j}(m)$ for some prime $p$ and integer $j$.
Recalling that $F_{k,j}(m)$ has no prime divisors smaller than $(\log x)^\epsilon$,
then by $(i)$ and $(ii)$ we must have $p>2MP/z^{1/2}$.
Hence, if we write $F_{k,j}(m)=\ell p^2$, it follows
\[
\ell = \frac{F_{k,j}(m)}{p^2}\leq  \frac{(mP)^2+j}{p^2}\leq \frac{2M^2P^2}{(2MP/z^{1/2})^2} <z.
\]
However, we have seen that $F_{k,j}(m)$ is not divisible by any prime $p \leq z$, and so $\ell=1$. In other words, we must have
\[
F_{k,j}(m) = p^2.
\]
This implies that 
\[
1 = p^2 - (mP)^2/j= p^2 - j(mP_j)^2
\]
a solution to Pell's equation with discriminant $j$.
The number of such $m\leq M$ is $\asymp_k \log M$.
After discarding these integers we obtain the lemma.
\end{proof}

\begin{proof}[Proof of Theorem \ref{TG}]
Let $m_0\pmod{q}$ be as in Lemma \ref{lemma1}.
Lemma~\ref{lemma2} provides $M^{1-o(1)}$ integers $m\leq M$ such that $m\equiv m_0\pmod{q}$ and
\[
D_j \gg_k \Delta(m)\qquad\forall\;j=1,\dots,k.
\]
Equivalently, since $\Delta(m)\asymp_k m^2$, this holds for $x^{1/2-o(1)}$ integers up to $x$.
Moreover, Lemma \ref{lemma1} shows that the condition $m\equiv m_0\pmod{q}$ ensures that
\[
\chi_{D_j}(p) = 1 \qquad \forall\; 5^k<p\leq (\log x)^\epsilon.
\]
Using this in \eqref{eq:L} (and recalling the remark that follows it),
we deduce
\begin{equation}\label{0912:eq001}
L(1,\chi_{D_j}) \gg \prod_{p\leq (\log x)^\epsilon}\left(1-\frac{\chi_{D_j}(p)}{p}\right)^{-1}\gg_k \log \log \Delta(m)
\end{equation}
with at most $O(x^{1/3})$ exceptions.
Finally, since by \eqref{bound:regulator} the regulator is bounded by $O(\log \Delta(m))$,
applying the class number formula we finish the proof:
\[
h_{\QQ(\sqrt{\Delta(m)+j})} = \frac{\sqrt{D_j}L(1,\chi_{D_j})}{\log\eps_{D_j}} \gg_k \frac{\sqrt{\Delta(m)}}{\log \Delta(m)} \log\log \Delta(m).
\]
\end{proof}

\section{Uniformity}\label{S3}

We can keep track of the dependence on $k$ throughout the proof of Theorem~\ref{TG}.
First, \eqref{bound:regulator} gives $\log \eps_{D_j}\ll \log \Delta(m)$ uniformly in $k$.

Moreover, we have $D_j\geq (s_jr_j)^{-1}\Delta(m) \geq s_j^{-1}\Delta(m)$
in the notation of Lemma \ref{lemma2}, where we used 
the definition of $D_j$ and the fact that $r_j=1$
if $m$ is selected as in Lemma \ref{lemma2}.
Noting that $s_j\leq \sqrt{k}$ and applying the class number formula, we deduce
\begin{equation}\label{0912:eq010}
h_{\QQ(\sqrt{\Delta(m)+j})} = \frac{\sqrt{D_j} L(1,\chi_{D_j})}{\log \eps_{D_j}}
\gg
\frac{\sqrt{\Delta(m)}L(1,\chi_{D_j})}{s_j \log \Delta(m)}
\gg
\frac{\sqrt{\Delta(m)}L(1,\chi_{D_j})}{\sqrt{k} \log \Delta(m)}.
\end{equation}

As for the $L$-value, using \eqref{0912:eq001}
and recalling that $\chi_{D_j}(p)=1$ for $5^k<p\leq (\log x)^\epsilon$ by Lemma \ref{lemma1}
(thus implicitly assuming $k\ll \log\log x$), we have
\begin{equation}\label{1812:eq001}
\begin{split}
L(1,\chi_{D_j})
&\gg
\exp\bigg(\sum_{p\leq (\log x)^\epsilon} \frac{\chi_{D_j}(p)}{p} + O(1)\bigg)
\\
&\gg
\exp\bigg(-\sum_{p\leq 5^k} \frac{1}{p} + \sum_{5^k<p\leq (\log x)^\epsilon} \frac{1}{p} + O(1)\bigg)
\gg
\frac{\log \log x}{k^2},
\end{split}
\end{equation}
where the implied constant is absolute.
Plugging the above into \eqref{0912:eq010} and recalling
that $x\geq \Delta(m)$ we obtain
\[
h_{\QQ(\sqrt{\Delta(m)+j})} \gg \frac{\sqrt{\Delta(m)}\log\log \Delta(m)}{k^{5/2}\log\Delta(m)},
\]
for every $j=1,\dots,k$, with an absolute implied  constant. 

\medskip

To improve the bound on average 
when we take the product of class numbers, we use again the first inequality in \eqref{0912:eq010}
and the estimate \eqref{1812:eq001} to write
\begin{equation}\label{1812:eq003}
\prod_{j=1}^k h_{\QQ(\sqrt{\Delta(m)+j})}
\gg
\bigg(\prod_{j=1}^k \frac{1}{s_j}\bigg) \bigg(\frac{\sqrt{\Delta(m)}\log\log \Delta(m)}{k^2 \log \Delta(m)}\bigg)^k.
\end{equation}
Regarding the product of $s_j$, we can estimate
\[
\begin{split}
\log \prod_{j=1}^k s_j = \sum_{j=1}^k \log s_j
&= \sum_{p^e} (\log p)\cdot |\{j\leq k:\; p^{2e}|j\}|
\\
&=
\sum_{p\leq k} (\log p)\left(\left\lfloor\frac{k}{p^2}\right\rfloor + \left\lfloor\frac{k}{p^4}\right\rfloor + \cdots\right)
\leq
k \sum_{p} \frac{\log p}{p^2-1} < k.
\end{split}
\]
In other words, $\prod_{j=1}^k s_j^{-1} > e^{-k}$.
Using this in \eqref{1812:eq003} we get
\[
\prod_{j=1}^k h_{\QQ(\sqrt{\Delta(m)+j})}
\gg
\frac{\left(\sqrt{\Delta(m)}\log\log \Delta(m)\right)^k}{e^{k} k^{2k} (\log \Delta(m))^k},
\]
which proves \eqref{1812:eq002}.

%% BIBLIOGRAPHY %%

\end{document}